\documentclass[12pt]{article}
\usepackage{amsbsy,amsfonts}

\newtheorem{theo}{Theorem}[section]
\newtheorem{lemma}[theo]{Lemma}
\newtheorem{prop}[theo]{Proposition}

\newtheorem{coll}[theo]{Corollary}
\newtheorem{remark}[theo]{Remark}

\newcommand{\be}{\begin{equation}}
\newcommand{\ee}{\end{equation}}
\newcommand{\ba}{\begin{array}}
\newcommand{\ea}{\end{array}}

\newcommand{\dsum}{\displaystyle \sum}

\begin{document}

\title{A conjecture of Herzog and Conca on counting of paths }
\author{ Hsin-Ju Wang \thanks{e-mail:  hjwang@math.ccu.edu.tw} \\
Department of Mathematics, National Chung Cheng University\\
          Chiayi 621, Taiwan}

\date{}

\maketitle

\begin{abstract}
~\\ A formula concerning counting of paths was conjectured by
Herzog and Conca few years ago. Recently, Krattenthaler and
Prohaska gave an affirmative answer to this conjecture. In this
paper we generalize this formula.
\end{abstract}
\section{Introduction}
Let $X$ be the set $\{(i, j)~|~1\leq i\leq m, 1\leq j\leq n\}$ of
the plane with the partial order given by $(i, j)\leq (i', j')$ if
$i\geq i'$ and $j\leq j'$. Let $P, Q\in X$ with $P\geq Q$; a path
from $P$ to $Q$ is a maximal chain in $X$ with end points $P$ and
$Q$. A corner of a path $C$ is an element $(i, j)\in C$ for which
$(i-1, j)$ and $(i, j-1)$ belong to $C$ as well. We use $w(P, Q)$
for the number of different paths from $P$ to $Q$ and $w_k(P, Q)$
for  the number of different paths with $k$ corners from $P$ to
$Q$. Therefore $w(P, Q)=\sum_{k\geq 0}w_k(P, Q)$.
\\Let $P_i$,$Q_i$, $i=1, \dots, r$ be points of $X$. A subset
$W\subseteq X$ is called an $r$-tuple of non-intersecting paths
form $P_i$ to $Q_i$ ($i=1, \dots, r$) if $W=C_1\cup C_2\cup \cdots
\cup C_r$ where each $C_i$ is a path from $P_i$ to $Q_i$, and
where $C_i\cap C_j=\emptyset$ if $i\neq j$. The number of corners
$c(W)$ of $W$ is the sum of the number of corners of the $C_i$. We
use $w_k({\bf P}, {\bf Q})$ for the number of the families of
non-intersecting paths form $P_i$ to $Q_i$ ($i=1, \dots, r$) with
exactly $k$ corners, and $W({\bf P}, {\bf Q})$ for the polynomial
(in $t$) $\sum_{k\geq 0} w_k({\bf P}, {\bf Q})t^k$. The work of
Krattenthaler \cite{kr} and Kulkarni \cite{ku2} showed the
following.
\begin{theo}
Let $X=\{(i, j)~|~1\leq i\leq m, 1\leq j\leq n\}$ with the partial
order giving as the above. Let $P_i=(a_i, n)$, $1= a_1< \cdots <
a_r\leq m$, and $Q_i=(m, b_i)$, $1=b_1< \cdots <b_r\leq n$. Then
$$W({\bf P}, {\bf Q}) = \det[\sum_{k\geq 0}{m-a_i+i-j\choose
k}{n-b_j+j-i\choose k+j-i}t^k]_{i,j=1, \dots, r}. \label{theo1}$$
\end{theo}
When $a_i=b_i=i$, the above formula is equivalent to the following
one obtained by Conca and Herzog \cite{cohe}, $$ W({\bf P}, {\bf
Q}) = t^{-{r\choose 2}} \det[\sum_{k\geq 0} {m-i\choose k
}{n-j\choose k}t^k]_{i,j=1, \dots, r}.$$ Since ${m-i\choose k
}{n-j\choose k}$ is the number of paths from $(i, n)$ to $(m, j)$
with exactly $k$ corners, this formula suggests that if $P_i=(i,
n)$ and $Q_i=(m, i)$ are points of a rectangular region $X$, then
\be \label{eq2} W({\bf P}, {\bf Q}) = t^{-{r\choose 2}}
\det[W(P_i, Q_j)]_{i,j=1, \dots, r}.\ee According to this formula,
Conca and Herzog made the following conjecture. \vskip 0.05in ~~\\
{\bf Conjecture:} Let $Y$ be a one-sided ladder-shaped region of
the plane (see Figure 1 below) with the partial order defined as
the above. Let $P_i=(i, n)$ and $Q_i=(m ,i)$ be points of $Y$.
Then $$W_Y({\bf P}, {\bf Q}) = t^{-{r\choose 2}}\det[W_Y(P_i,
Q_j)]_{i,j=1, \dots, r}.$$ \vskip 0.08in
\begin{picture}(0,0)

\put(120,0){\line(0,-1){100}}\put(120,-100){\line(1,0){140}}
\put(120,0){\line(1,0){80}}\put(200,0){\line(0,-1){20}}
\put(200,-20){\line(1,0){40}}\put(240,-20){\line(0,-1){30}}
\put(240,-50){\line(1,0){20}}\put(260,-50){\line(0,-1){50}}

\put(115,-5){$\times$} \put(105,5){$P_1$} \put(125,-5){$\times$}
\put(118,5){$P_2$} \put(165,-5){$\times$} \put(155,5){$P_r$}

\put(257,-105){$\times$} \put(262,-115){$Q_1$}

\put(257,-60){$\times$} \put(262,-70){$Q_r$}

\put(150,-125){Figure~1}

\end{picture}
\vskip 1.8in ~~\\ Recently, Krattenthaler and Prohaska \cite{krpr}
gave an affirmative answer to the Conjecture by using the notion
{\it two-rowed arrays} introduced in \cite{krmo}. \par Our main
concern in this paper is trying to give a self-contained proof to
the result obtained by  Krattenthaler and Prohaska. As a
consequence, we are able to generalize their result as follows.
\begin{theo}
Let $Y$ be a one-sided ladder-shaped region of the plane with the
partial order given in the beginning. Let $P_i=(i, n)$ and $Q=(m,
b_i)$, $i=1, \dots, r$ be points of $Y$. Let $W(n,m; b_1, \dots,
b_r)$ be the polynomial $\sum_{k\geq 0} w_k({\bf P}, {\bf Q})t^k$.
Let $\tilde W(n,m; b_1, \dots, b_r)=\det[W(P_i, Q_j)]_{i,j=1,
\dots, r}$. Then $$t^{-{r\choose 2}}\tilde W(n,m; b_1, \dots,
b_r)= \sum_{\underline{c}\in S_r} A(\underline{b};
\underline{c})(1-t)^{\sum (c_i-b_i)}W(n,m; c_1, \dots, c_r),$$
where $\underline{b}=(b_1, \dots, b_r)$ and $\underline{c}=(c_1,
\dots, c_r)$. In particular, if $b_i=i~\forall i$, then
$$t^{-{r\choose 2}}\tilde W(n,m; 1, \dots, r)= W(n,m; 1, \dots,
r).$$ \label{theo2}
\end{theo}
Here $S_r=\{(b_1, b_2, \dots , b_r)\in \mathbb{N}^r ~|~b_1<b_2<
\cdots <b_r\}$ and $A(\underline{b}; \underline{c})$ is defined in
section 3.
\section{Some fundamental lemmas}
Let $Y$ be a one-sided ladder-shaped region of the plane with the
partial order given by $(i, j)\leq (i', j')$ if $i\geq i'$ and
$j\leq j'$. Let $P_i=(i, s)$ and $Q_i=(l, b_i)$, $i=1, \dots, r$
be points of $Y$. We use $W(s,l;b_1, \dots, b_r)$ for the
polynomial $\sum_{k\geq 0} w_k({\bf P}, {\bf Q})t^k$ and $\tilde
W(s,l;b_1, \dots, b_r)$ for $\det[W(P_i, Q_j)]_{i,j=1, \dots, r}$.
In order to obtain some useful properties of $W$ and $\tilde W$,
we introduce the following notations. \\ Let $l$ be a positive
integer and $$S_l=\{(b_1, b_2, \dots , b_l)\in \mathbb{N}^l
~|~b_1<b_2< \cdots <b_l\}.$$ We define a partial order on $S_l$ so
that for any $(b_1, \dots, b_l), (c_1, \dots, c_l)\in S_l$, $(b_1,
\dots, b_l)\leq (c_1, \dots, c_l)$ if $b_i\leq c_i~\forall i$.
Moreover, for any two lattice points $(b_1, \dots, b_l), (c_1,
\dots, c_l)\in S_l$, we define $\{ \ba{ccc} c_1 & \cdots & c_l \\
b_1 & \cdots & b_l \ea \}$ and $[ \ba{ccc} c_1 & \cdots & c_l \\
b_1 & \cdots & b_l \ea ]$ as follows. Let $b, c\in \mathbb{N}$. We
use the symbol $[\ba{c} c \\ b \ea ]$ for $t$, $1$ and $0$ if
$c>b$, $c=b$ and $c<b$, where $t$ is a variable. Furthermore, for
any points $(b_1, \dots, b_r), (c_1, \dots, c_r)\in S_r$,
$$[\ba{ccc} c_1 & \cdots & c_r \\ b_1 & \cdots & b_r \ea ]=0$$ if
$c_i\geq b_{i+1}$ for some $i<r$ and $$[\ba{ccc} c_1 & \cdots &
c_r \\ b_1 & \cdots & b_r \ea ]=[\ba{c} c_1 \\ b_1 \ea ] \cdots
[\ba{c} c_r
\\ b_r \ea ]$$ if $c_i<b_{i+1}$ for every $i$. Finally,
$$\{\ba{ccc} c_1 & \cdots & c_r \\ b_1 & \cdots & b_r \ea \}=
\sum_{\sigma\in S_X} sign(\sigma)[\ba{c} \sigma(c_1) \\ b_1 \ea
]\cdots [\ba{c} \sigma(c_r) \\ b_1 \ea ],$$ where $X=\{c_1, \dots,
c_r)$ and $S_X$ is the permutation group on $X$.
\begin{lemma}
Let $Y$ be a one-sided ladder-shaped region of the plane with the
partial order given  in the beginning. Suppose that $m=r$ (Then
$Y$ is a rectangular region). Let $P_i=(i, n)$ and $Q_i=(r, b_i)$,
$i=1, \dots, r$. Then $$t^{-{r\choose 2}}\tilde W(n,r; b_1, \dots,
b_r)=w({\bf P}, {\bf Q}).$$ \label{lem1}
\end{lemma}
\begin{proof}
Notice that $$\ba{rcl} W(P_i, Q_j) & = & \sum_{k\geq 0}
{r-i\choose k}{n-b_j\choose k}t^k \\ & = &  \sum_{k\geq 0}^{r-i-1}
{r-i\choose k}{n-b_j\choose k}t^k+ {n-b_j\choose r-i}t^{r-i}, \ea
$$ therefore, the $i$-th row of $[W(P_i, Q_j)]_{i,j=1, \dots, r}$
is $$\ba{rcl} & & \\ <\sum_{k\geq 0}^{r-i-1} {r-i\choose
k}{n-b_1\choose k}t^k, \dots, \sum_{k\geq 0}^{r-i-1} {r-i\choose
k}{n-b_r\choose k}t^k>+<{n-b_1\choose r-i}t^{r-i}, \dots,
{n-b_r\choose r-i}t^{r-i}>.\ea $$ Since $$\ba{rcl} & & \\
<\sum_{k\geq 0}^{r-i-1} {r-i\choose k}{n-b_1\choose k}t^k, \dots,
\sum_{k\geq 0}^{r-i-1} {r-i\choose k}{n-b_r\choose k}t^k>\ea $$ is
a linear combination of the last $r-i$ rows of $[W(P_i,
Q_j)]_{i,j=1, \dots, r}$, one can use elementary row operations to
obtain that $$\ba{rcl} \tilde W(n,r; b_1, \dots, b_r) &
 = & \det \left( \ba{ccc} W(P_1, Q_1) & \cdots &  W(P_1, Q_r) \\
\cdots & \cdots & \cdots \\ 1+(n-b_1)t & \cdots & 1+(n-b_r)t \\ 1
& \cdots & 1 \ea  \right) \\ & = & \det \left( \ba{ccc} W(P_1,
Q_1) & \cdots &  W(P_1, Q_r) \\ \cdots & \cdots & \cdots \\
(n-b_1)t & \cdots & (n-b_r)t \\ 1 & \cdots & 1 \ea  \right) \\ & =
& \cdots  \\  & = & \det \left( \ba{ccc} W(P_1, Q_1) & \cdots &
W(P_1, Q_r) \\ \cdots & \cdots & \cdots \\  \sum_{k\geq 0}
{r-i\choose k}{n-b_1\choose k}t^k  & \cdots & \sum_{k\geq 0}
{r-i\choose k}{n-b_r\choose k}t^k
\\ {n-b_1\choose r-i-1}t^{r-i-1}  & \cdots &
{n-b_r\choose r-i-1}t^{r-i-1} \\ \cdots & \cdots & \cdots \\ 1 &
\cdots & 1  \ea \right) \\ & = &  \det \left( \ba{ccc} W(P_1, Q_1)
& \cdots & W(P_1, Q_r) \\ \cdots & \cdots & \cdots \\
{n-b_1\choose r-i}t^{r-i}  & \cdots & {n-b_r\choose r-i}t^{r-i}
\\ {n-b_1\choose r-i-1}t^{r-i-1}  & \cdots &
{n-b_r\choose r-i-1}t^{r-i-1} \\ \cdots & \cdots & \cdots \\ 1 &
\cdots & 1  \ea \right) \\ & = & \cdots \\ & = & \det \left(
\ba{ccc} {n-b_1\choose r-1}t^{r-1} & \cdots & {n-b_r\choose
r-1}t^{r-1} \\ \cdots & \cdots & \cdots \\ 1 & \cdots & 1  \ea
\right) \\ & = & ct^{{r\choose 2}},   \ea $$ for some constant
$c$. However, $\tilde W(n,r; b_1, \dots, b_r)(1)=w({\bf P}, {\bf
Q})$, therefore the assertion follows.
\end{proof}
\begin{coll}
Let $Y$ be a one-sided ladder-shaped region of the plane with the
partial order given  in the beginning. Suppose that $m=r$ (Then
$Y$ is a rectangular region). Let $P_i=(i, n)$ and $Q_i=(r, b_i)$,
$i=1, \dots, r$. Then $$t^{-{r\choose 2}}\tilde W(n,r; b_1, \dots,
b_r)=t^{-{r-1\choose 2}}\sum_{c_1=b_1}^{b_2-1} \cdots
\sum_{c_{r-1}=b_{r-1}}^{b_r-1} \tilde W(b_r-1,r-1; c_1, \dots,
c_{r-1}).$$ \label{coll1}
\end{coll}
\begin{proof}
Let $P'_i=(i, b_r-1)$ and $S_{c_i}=(r-1, c_i)$. Then by
Lemma~\ref{lem1}, it is enough to show that $$w({\bf P}, {\bf Q})=
\sum_{c_1=b_1}^{b_2-1} \cdots \sum_{c_{r-1}=b_{r-1}}^{b_r-1}w({\bf
P'}, {\bf S_{\underline{c}}}), $$ where ${\bf P'}=\{P'_1, \dots,
P'_{r-1}\}$ and ${\bf S_{\underline{c}}}=\{S_{c_1}, \dots,
S_{c_{r-1}}\}$. However, this follows from two facts. One is $$
w({\bf P}, {\bf Q})= w({\bf P'}, {\bf Q'}),$$ where ${\bf
Q'}=\{Q_1, \dots, Q_{r-1}\}$. The other is $$\ba{rcl}  w({\bf P'},
{\bf Q'}) & = & \sum_{c_1=b_1}^{b_2-1} w({\bf P'}, S_{c_1},Q_2,
\dots, Q_{r-1})\\ & = & \cdots \\ & = & \sum_{c_1=b_1}^{b_2-1}
\cdots \sum_{c_{r-1}=b_{r-1}}^{b_r-1}w({\bf P'}, {\bf
S_{\underline{c}}}) \ea $$ by \cite[Lemma~3.2]{wa2}
\end{proof}
\begin{lemma}
Let $Y$ be a one-sided ladder-shaped region of the plane with the
partial order given  in the beginning. Then the following hold:
\begin{description}
\item{(i)} $$\tilde W(n,m;b_1, \cdots,
b_r)=\sum_{\underline{d}\in S_r,(m,d_r)\in Y} \{ \ba{ccc} d_1 &
\cdots & d_r \\ b_1 & \cdots & b_r \ea \} \tilde W(n,m-1;d_1,
\dots, d_r).$$
\item{(ii))}
$$W(n,m;b_1, \cdots, b_r)=\sum_{\underline{d}\in S_r,(m,d_r)\in Y}
[ \ba{ccc} d_1 & \cdots & d_r \\ b_1 & \cdots & b_r \ea ]
W(n,m-1;d_1, \dots, d_r).$$
\end{description} \label{lem2}
\end{lemma}
\begin{proof}
(i). Let $S_{d_j}=(m-1, d_j)$. Since $$W(P_i, Q_j)=\sum_{d_j\in
\mathbb{N},(m,d_j)\in Y}[\ba{c} d_j \\ b_j \ea ] W(P_i, S_{d_j})$$
by \cite[Lemma~3.1]{wa2}, $$\ba{rcl} & &  \tilde W(n,m;b_1,
\cdots, b_r)\\ & = & \det[W(P_i, Q_j)]_{i,j=1, \dots, r } \\ & = &
\dsum_{d_1\in \mathbb{N},(m,d_1)\in Y}\cdots \sum_{d_r\in
\mathbb{N},(m,d_r)\in Y} \det \{ [\ba{c} d_j \\ b_j \ea ] W(P_i,
S_{d_j}) \}_{i,j=1, \dots, r} \\ & = & \dsum_{d_1\in
\mathbb{N},(m,d_1)\in Y}\cdots \sum_{d_r\in \mathbb{N},(m,d_r)\in
Y} [\ba{c} d_1 \\ b_1 \ea ]\cdots [\ba{c} d_r \\ b_r \ea ] \det [
W(P_i, S_{d_j}) ]_{i,j=1, \dots, r}\\ & = &
\dsum_{\underline{d}\in S_r,(m,d_r)\in Y} \sum_{\sigma}[\ba{c}
\sigma(d_1)
\\ b_1 \ea ]\cdots [\ba{c} \sigma(d_r) \\ b_r \ea ] \det [ W(P_i,
S_{\sigma(d_j)}) ]_{i,j=1, \dots, r} \\  & = &
\dsum_{\underline{d}\in S_r,(m,d_r)\in Y} \sum_{\sigma}
sign(\sigma) [\ba{c} \sigma(d_1)
\\ b_1 \ea ]\cdots [\ba{c} \sigma(d_r) \\ b_r \ea ] \det [ W(P_i,
S_{d_j}) ]_{i,j=1, \dots, r}\\ & = & \dsum_{\underline{d}\in
S_r,(m,d_r)\in Y} \{ \ba{ccc} d_1 & \cdots & d_r \\ b_1 & \cdots &
b_r \ea \} \tilde W(n,m-1;d_1, \dots, d_r). \ea $$ (ii). Let
$S_{d_j}=(m-1, d_j)$. By using \cite[Lemma~3.2]{wa2} repeatedly,
one can obtain that $$W(n,m;b_1, \cdots,
b_r)=\dsum_{d_1=b_1}^{b_2-1}\cdots \sum_{d_r=b_r}^b W(n,m-1;d_1,
\dots,d_r),$$ where $b$ is the maximal integer for which $(m,
b)\in Y$. Thus (ii) follows from the definition of $[ \ba{ccc} d_1
& \cdots & d_r \\ b_1 & \cdots & b_r \ea ]$.
\end{proof}

\section{Properties of $A(\underline{b}, \underline{c})$}
Let $l\geq 2$ be an integer and $$S_l=\{(b_1, b_2, \dots , b_l)\in
\mathbb{N}^l ~|~b_1<b_2< \cdots <b_l\}$$ with the partial order
given in section 2.\\ For any two lattice points $(b_1, \dots,
b_l), (c_1, \dots, c_l)\in S_l$, we define $A(b_1, \dots, b_l;
c_1, \dots, c_l)$ as follows.
\\ Assume $l=2$. Then $A(b_1,b_2;c_1,c_2)=0$ if $b_2\neq c_2$ or
$c_1<b_1$, and $A(b_1,b_2;c_1,c_2)=1$ if $b_2=c_2$ and $b_1\leq
c_1$.
\\ Assume that $l>2$. If $b_l\neq c_l$,
then $ A(b_1, \dots, b_l; c_1, \dots,c_l)=0$. If $b_l=c_l$, then
$$A(b_1, \dots, b_l; c_1, \dots,c_l)=\sum_{\underline{c'}\leq
\underline{c},\underline{c'}\in S_{l-1}} A(b_1, \dots, b_{l-1};
c'_1, \dots,c'_{l-1}).$$
\begin{theo}
Let $l\geq 2$ be a positive integer and $(b_1, \dots, b_l), (c_1,
\dots, c_l)\in S_l$. Let $$B_1(b_1, \dots, b_l;c_1, \dots,
c_l)=\sum_{\underline{d}\in S_l} A(d_1, \dots, d_l;c_1, \dots,
c_l) \{ \ba{rcl} d_1 & \cdots & d_l
\\ b_1 & \cdots & b_l \ea \}(1-t)^{\sum (c_i-d_i)}, $$ $$B_2(b_1,
\dots, b_l;c_1, \dots, c_l)=\sum_{\underline{d}\in S_l} A(b_1,
\dots, b_l;d_1, \dots, d_l) [ \ba{rcl} c_1 & \cdots & c_l
\\ d_1 & \cdots & d_l \ea ](1-t)^{\sum (d_i-b_i)} $$ and
$$B_3(b_1, \dots, b_l;c_1, \dots, c_l)=\dsum_{d_1=b_1}^{b_2-1}
\cdots \dsum_{d_{l-1}=b_{l-1}}^{b_l-1}A(d_1, \dots, d_{l-1};c_1,
\dots, c_{l-1})(1-t)^{\sum_{i=1}^{l-1} (c_i-d_i)}$$ if $l\geq 3$
and $b_l=c_l$. Then $$B_1(b_1, \dots, b_l;c_1, \dots,
c_l)=B_2(b_1, \dots, b_l;c_1, \dots, c_l)$$ and $$B_1(b_1, \dots,
b_l;c_1, \dots, c_l)=B_3(b_1, \dots, b_l;c_1, \dots, c_l)$$ if
$l\geq 3$ and $b_l=c_l$. \label{theo3}
\end{theo}

\begin{remark}
(i). $B_2(b_1, \dots, b_l;c_1, \dots, c_l)=0$ if $c_{l-1}\geq
b_l$.\\ (ii). $B_3(b_1, \dots, b_l;c_1, \dots, c_l)=B_3(b_1,
\dots, b_{l-1}, b;c_1, \dots, c_{l-1}, b)$ for every $b>
max\{b_{l-1}, c_{l-1}\}$. \label{rem1}
\end{remark}
To show Theorem~\ref{theo3}, we need several lemmas. The first one
is easy to obtain, we left the proof to the reader.

\begin{lemma}
Let $b\leq e\leq c$ be positive integers. Then $$\sum_{d=e}^c
[\ba{c} c \\ d \ea ] (1-t)^{d-b} =(1-t)^{e-b}$$ and $$\sum_{d=b}^e
[\ba{c} d \\ b \ea ] (1-t)^{c-d} =(1-t)^{c-e}.$$ \label{lem3}
\end{lemma}
\begin{lemma}
Let $l\geq 2$ be a positive integer. Let $(b_1, \dots, b_l)$ and
$(c_1, \dots, c_l)$ be two points in $S_l$. Then
$$\sum_{\underline{b}\leq \underline{d},\underline{d}\in S_l}
[\ba{rcl} c_ 1 & \dots & c_l \\ d_1 & \cdots & d_l \ea ]
(1-t)^{\sum (d_i-b_i)}= \sum_{\underline{d}\leq
\underline{c},\underline{d}\in S_l} [ \ba{rcl} d_1 & \dots & d_l
\\ b_1 & \cdots & b_l \ea ] (1-t)^{\sum (c_i-d_i)}.$$ \label{lem4}
\end{lemma}
\begin{proof}
We may assume that $(b_1, \dots, b_l)\leq (c_1, \dots, c_l)$. By
Lemma~\ref{lem3}, $$\ba{rcl} &  & \dsum_{\underline{b}\leq
\underline{d},\underline{d}\in S_l} [\ba{rcl} c_ 1 & \dots & c_l
\\ d_1 & \cdots & d_l \ea ] (1-t)^{\sum (d_i-b_i)} \\ & = &
\{\dsum_{d_1=b_1}^{c_1}[\ba{c} c_1 \\ d_1\ea ]
(1-t)^{d_1-b_1}\}\prod_{i=2}^l\{\sum_{d_i=max\{b_i,
c_{i-1}+1\}}^{c_i}[ \ba{c} c_i \\ d_i\ea ] (1-t)^{d_i-b_i}\} \\ &
= & \prod_{i=2}^l (1-t)^{max\{0, c_{i-1}-b_i+1\}} . \ea $$ Also,
by Lemma~\ref{lem3}, $$\ba{rcl} &  & \dsum_{\underline{d}\leq
\underline{c},\underline{d}\in S_l} [ \ba{rcl} d_1 & \dots & d_l
\\ b_1 & \cdots & b_l \ea ] (1-t)^{\sum (c_i-d_i)} \\ & = &
\prod_{i=1}^{l-1} \{\sum_{d_i=b_i}^{min\{c_i, b_{i+1}-1\}}[\ba{c}
d_i\\ b_i\ea ](1-t)^{c_i-d_i} \}\{\dsum_{d_l=b_l}^{c_l} [\ba{c}
d_l\\ b_l\ea ](1-t)^{c_l-d_l} \} \\ & = & \prod_{i=1}^{l-1}
(1-t)^{max\{0, c_i-b_{i+1}+1\}}. \ea $$ This shows the Lemma.
\end{proof}
\begin{lemma}
Let $(b_1,b_2), (c_1,c_2)\in S_2$. If $c_1\leq b_2-1$ and $b_1\leq
c_1$, then  $$B_1(b_1, b_2;c_1, c_2)=B_2(b_1,b_2;c_1,c_2)=[\ba{c}
c_2 \\ b_2 \ea ] .$$ If $c_1<b_1$ or $c_1\geq b_2$, then
$$B_1(b_1, b_2;c_1, c_2)=B_2(b_1,b_2;c_1,c_2)=0.$$ \label{lem5}
\end{lemma}
\begin{proof}
It is clear that if $c_1<b_1$ then $B_1(b_1, b_2;c_1,
c_2)=B_2(b_1,b_2;c_1,c_2)=0$, so that we may assume $b_1\leq c_1$.
If $c_1\leq b_2-1$, then by Lemma~\ref{lem3} $$\ba {rcl} B_1(b_1,
b_2;c_1, c_2) & = & \sum_{d_1=b_1}^{c_1} A(d_1,c_2;c_1,c_2) \{
\ba{cc} d_1 & c_2
\\ b_1 & b_2 \ea \} (1-t)^{c_1-d_1}\\ & = & \sum_{d_1=b_1}^{c_1}
[ \ba{c} d_1 \\ b_1 \ea ][ \ba{c} c_2 \\ b_2 \ea ](1-t)^{c_1-d_1}
 \\ & = & [ \ba{c} c_2 \\ b_2 \ea ], \ea $$ and $$\ba {rcl}
B_2(b_1, b_2;c_1, c_2) & = &
\sum_{d_1=b_1}^{c_1}A(b_1,b_2;d_1,b_2) [ \ba{cc} c_1 & c_2
\\ d_1 & b_2 \ea ] (1-t)^{d_1-b_1}\\ & = & \sum_{d_1=b_1}^{c_1}
[ \ba{c} c_1 \\ d_1 \ea ][ \ba{c} c_2 \\ b_2 \ea ](1-t)^{d_1-b_1}
\\ & = & [ \ba{c} c_2 \\ b_2 \ea ]. \ea $$ On the other hand, if
$c_1\geq b_2$ then $$\ba {rcl} B_1(b_1, b_2;c_1, c_2) & = &
\sum_{d_1=b_1}^{c_1}A(d_1,c_2;c_1,c_2) \{ \ba{cc} d_1 & c_2
\\ b_1 & b_2 \ea \} (1-t)^{c_1-d_1}\\ & = & \sum_{d_1=b_1}^{b_2-1}
[ \ba{c} d_1 \\ b_1 \ea ][ \ba{c} c_2 \\ b_2 \ea ](1-t)^{c_1-d_1}
+  \{ \ba{cc} b_2 & c_2
\\ b_1 & b_2 \ea \} (1-t)^{c_1-b_2} \\ & = & \sum_{d_1=b_1}^{b_2-1}
[ \ba{c} d_1 \\ b_1 \ea ]t(1-t)^{c_1-d_1} + (t^2-t)(1-t)^{c_1-b_2}
\\ & = & 0, \ea $$ and $$B_2(b_1, b_2;c_1, c_2)=
\sum_{d_1=b_1}^{b_2-1}A(b_1,b_2;d_1,b_2) [ \ba{cc} c_1 & c_2 \\
d_1 & b_2 \ea ] (1-t)^{d_1-b_1}= 0  $$ as $[ \ba{cc} c_1 & c_2 \\
d_1 & b_2 \ea ] =0$. This proves the lemma.
\end{proof}
\begin{lemma}
Let $(b_1, \dots, b_l)\in S_l$ with $l\geq 3$. Let $b\in
\mathbb{N}$ such that $b>b_l$. For every $i\leq l$, let $\sigma_i:
\{b_1, \dots, b_l\} \longrightarrow \{b_1, \dots, \hat{b}_i,
\dots, b_l, b\}$ be the bijective map such that
$\sigma_i(b_j)=b_j$ if $j<i$, $\sigma_i(b_j)=b_{j+1}$ if $i\leq
j\leq {l-1}$ and $\sigma_i(b_l)=b$. Let $f(x_1, \dots, x_{l-1})\in
\mathbb{Q}(t)[x_1, \dots, x_{l-1}]$. Then $$\ba{rcl} & &
\dsum_{i=1}^l (-1)^i \{
\sum_{d_1=\sigma_i(b_1)}^{\sigma_i(b_2)-1}\cdots \sum_{d_{l-1}=
\sigma_i(b_{l-1})}^{\sigma_i(b_l)-1} f(d_1, \dots, d_{l-1})\} \\ &
= & (-1)^l \dsum_{d_1=b_1}^{b_2-1}\cdots
\sum_{d_{l-1}=b_{l-1}}^{b_l-1} f(d_1, \dots, d_{l-1}).\ea $$
\label{lem6}
\end{lemma}
\begin{proof}
We proceed by induction on $l$. The case $l=3$ is easy to check.
We may assume that $l>3$. Then by induction $$ \ba{rcl}
 & & \dsum_{i=1}^l (-1)^i \{
\sum_{d_1=\sigma_i(b_1)}^{\sigma_i(b_2)-1}\cdots \sum_{d_{l-1}=
\sigma_i(b_{l-1})}^{\sigma_i(b_l)-1} f(d_1, \dots, d_{l-1})\} \\ &
= & (-1)^l \dsum_{d_1=b_1}^{b_2-1}\cdots
\sum_{d_{l-2}=b_{l-2}}^{b_{l-1}-1} \sum_{d_{l-1}=b_{l-1}}^{b-1}
f(d_1, \dots, d_{l-1}) \\ & & + \dsum_{i=1}^{l-1} (-1)^i
\sum_{d_{l-1}=b_l}^{b-1}\{
\sum_{d_1=\sigma_i(b_1)}^{\sigma_i(b_2)-1}\cdots \sum_{d_{l-2}=
\sigma_i(b_{l-2})}^{\sigma_i(b_{l-1})-1} f(d_1, \dots, d_{l-1})\}
\\ & = & (-1)^l \dsum_{d_1=b_1}^{b_2-1}\cdots
\sum_{d_{l-2}=b_{l-2}}^{b_{l-1}-1} \sum_{d_{l-1}=b_{l-1}}^{b-1}
f(d_1, \dots, d_{l-1}) \\ & & +
(-1)^{l-1}\dsum_{d_1=b_1}^{b_2-1}\cdots
\sum_{d_{l-2}=b_{l-2}}^{b_{l-1}-1} \sum_{d_{l-1}=b_l}^{b-1} f(d_1,
\dots, d_{l-1}) \\ & = & (-1)^l \dsum_{d_1=b_1}^{b_2-1}\cdots
\sum_{d_{l-1}=b_{l-1}}^{b_l-1} f(d_1, \dots, d_{l-1}).\ea $$
\end{proof}
 ~~\\ {\bf Proof of Theorem~3.1}
We prove the theorem by induction on $l$. If $l=2$, then it is the
content of Lemma~\ref{lem5}. Therefore we may assume that $l\geq
3$. Without loss of generality, we may further assume that $(b_1,
\dots, b_l)\leq (c_1, \dots, c_l)$. Let $$D(b_1, \dots,
b_{l-1};c_1, \dots, c_{l-1})=\sum_{\underline{d}\leq
\underline{c},\underline{d}\in S_{l-1}}B_1(b_1, \dots,
b_{l-1};d_1, \dots, d_{l-1})(1-t)^{\sum (c_i-d_i)}.$$ We will show
in the following that \be D(b_1, \dots, b_{l-1};c_1, \dots,
c_{l-1})=B_i(b_1, \dots, b_{l-1},b;c_1, \dots, c_{l-1},b)
\label{eq2} \ee for $i=1, 2, 3$ and for every $b> max\{b_{l-1},
c_{l-1}\}$.
\\ Observe first that $$ \ba{rcl} & & B_1(b_1,
\dots, b_{l-1},b;c_1, \dots, c_{l-1},b)\\ & = &
\dsum_{\underline{d}\in S_{l-1}} A(d_1, \dots, d_{l-1}, b;c_1,
\dots, c_{l-1}, b) \{ \ba{rcl} d_1 & \cdots & d_{l-1}
\\ b_1 & \cdots & b_{l-1} \ea \}(1-t)^{\sum_{i=1}^{l-1} (c_i-d_i)}\\ & = &
\dsum_{\underline{d}\in S_{l-1}} \sum_{\underline{d'}\leq
\underline{c}} A(d_1, \dots, d_{l-1};d'_1, \dots, d'_{l-1}) \{
\ba{rcl} d_1 & \cdots & d_{l-1}
\\ b_1 & \cdots & b_{l-1} \ea \}(1-t)^{\sum (c_i-d_i)}\\ & = &
\dsum_{\underline{d'}\leq \underline{c}}[\sum_{\underline{d}\in
S_{l-1}} A(d_1, \dots, d_{l-1};d'_1, \dots, d'_{l-1}) \{ \ba{rcl}
d_1 & \cdots & d_{l-1}
\\ b_1 & \cdots & b_{l-1} \ea \}(1-t)^{\sum (d'_i-d_i)}]
(1-t)^{\sum (c_i-d'_i)}\\& = &  \dsum_{\underline{d'}\leq
\underline{c}}B_1(b_1, \dots, b_{l-1};d'_1, \dots,
d'_{l-1})(1-t)^{\sum (c_i-d'_i)}\\ & = & D(b_1, \dots,
b_{l-1};c_1, \dots, c_{l-1}). \ea $$ Moreover, by induction and
Lemma~\ref{lem4} $$\ba{rcl} & & B_2(b_1, \dots, b_{l-1},b;c_1,
\dots, c_{l-1},b) \\ & = & \dsum_{\underline{d}\in S_{l-1}} A(b_1,
\dots, b_{l-1},b;d_1, \dots, d_{l-1}, b) [ \ba{cccc} c_1 & \cdots
& c_{l-1} & b \\ d_1 & \cdots & d_{l-1} & b \ea ]
(1-t)^{\sum_{i=1}^{l-1} (d_i-b_i)} \\ & = &
\dsum_{\underline{d}\in S_{l-1}} \sum_{\underline{d'}\leq
\underline{d}}A(b_1, \dots, b_{l-1};d'_1, \dots, d'_{l-1})[
\ba{rcl} c_1 & \cdots & c_{l-1} \\ d_1 & \cdots & d_{l-1}\ea ]
(1-t)^{\sum (d_i-b_i)} \\ & = & \dsum_{\underline{d'}\in S_{l-1}}
A(b_1, \dots, b_{l-1};d'_1, \dots, d'_{l-1})
\{\sum_{\underline{d'}\leq \underline{d}}[ \ba{rcl} c_1 & \cdots &
c_{l-1} \\ d_1 & \cdots & d_{l-1}   \ea ](1-t)^{\sum (d_i-d'_i)}
\} (1-t)^{\sum (d'_i-b_i)}
\\ & = & \dsum_{\underline{d'}\in S_{l-1}} A(b_1, \dots,
b_{l-1};d'_1, \dots, d'_{l-1})\{\sum_{\underline{d}\leq
\underline{c}}[ \ba{rcl} d_1 & \cdots & d_{l-1} \\ d'_1 & \cdots &
d'_{l-1}   \ea ](1-t)^{\sum (c_i-d_i}) \} (1-t)^{\sum (d'_i-b_i)}
\\ & = & \dsum_{\underline{d}\leq \underline{c}} \{
\sum_{\underline{d'}\in S_{l-1} }A(b_1, \dots, b_{l-1};d'_1,
\dots, d'_{l-1})[ \ba{rcl} d_1 & \cdots & d_{l-1} \\ d'_1 & \cdots
& d'_{l-1}   \ea ](1-t)^{\sum (d'_i-b_i)} \}(1-t)^{\sum (c_i-d_i)}
\\ & = & \dsum_{\underline{d}\leq \underline{c},\underline{d}\in
S_{l-1}}B_2(b_1, \dots, b_{l-1};d_1, \dots,
d_{l-1})(1-t)^{\sum_{i=1}^{d-1} (c_i-d_i)}.\\ & = & D(b_1, \dots,
b_{l-1};c_1, \dots, c_{l-1}). \ea $$ Finally, by induction and
Lemma~\ref{lem3} $$\ba{rcl}
 & & D(b_1, \dots, b_{l-1};c_1, \dots, c_{l-1}) \\ & = &
\dsum_{\underline{d}\leq \underline{c},\underline{d}\in
S_{l-1}}B_2(b_1, \dots, b_{l-1};d_1, \dots,
d_{l-1})(1-t)^{\sum_{i=1}^{l-1} (c_i-d_i)}\\ & = &
\dsum_{\underline{d}\leq \underline{c}} \{ \sum_{\underline{d'}\in
S_{l-1} }A(b_1, \dots, b_{l-1};d'_1, \dots, d'_{l-1})[ \ba{rcl}
d_1 & \cdots & d_{l-1} \\ d'_1 & \cdots & d'_{l-1}   \ea
](1-t)^{\sum (d'_i-b_i)} \}(1-t)^{\sum (c_i-d_i) } \\ & = &
\dsum_{\underline{d}\leq \underline{c}} \{ \sum_{\underline{d'}\in
S_{l-2} }A(b_1, \dots, b_{l-1};d'_1, \dots, d'_{l-2}, b_{l-1})[
\ba{rcl} d_1 & \cdots & d_{l-1} \\ d'_1 & \cdots & b_{l-1}   \ea
](1-t)^{\sum (d'_i-b_i)} \}(1-t)^{\sum (c_i-d_i)} \\ & = &
\dsum_{\underline{d}\leq \underline{c}, d_{l-2}< b_{l-1}} \{
\sum_{\underline{d'}\in S_{l-2}} A(b_1, \dots, b_{l-1};d'_1,
\dots, d'_{l-2}, b_{l-1})[ \ba{rcl} d_1 & \cdots & d_{l-2}
\\ d'_1 & \cdots & d'_{l-2}   \ea ](1-t)^{\sum (d'_i-b_i)} \} \\  &
& \cdot \{\sum_{d_{l-1}=b_{l-1}}^{c_{l-1}}[\ba{c} d_{l-1} \\
b_{l-1} \ea ](1-t)^{c_{l-1}-d_{l-1}} \} \cdot
(1-t)^{\sum_{i=1}^{l-2} (c_i-d_i)} \\ & = &
\dsum_{\underline{d}\leq \underline{c}}B_2(b_1, \dots,
b_{l-1};d_1, \dots, d_{l-2}, b_{l-1})(1-t)^{\sum_{i=1}^{l-2}
(c_i-d_i)}. \ea $$ If $l=3$, then by Lemma~\ref{lem5} $$\ba{rcl} &
& \dsum_{\underline{d}\leq \underline{c}} B_2(b_1, \dots,
b_{l-1};d_1, \dots, d_{l-2}, b_{l-1})(1-t)^{\sum_{i=1}^{l-2}
(c_i-d_i)} \\ & = & \dsum_{\underline{d}\leq \underline{c}}
B_2(b_1,b_2;d_1,b_2)(1-t)^{c_1-d_1} \\ & = &
\dsum_{d_1=b_1}^{min\{c_1,b_2-1 \}}
B_2(b_1,b_2;d_1,b_2)(1-t)^{c_1-d_1} \\ & = &
\dsum_{d_1=b_1}^{min\{c_1,b_2-1 \}} (1-t)^{c_1-d_1} \\  & = &
\dsum_{d_1=b_1}^{b_2-1} A(d_1,c_2;c_1,c_2)(1-t)^{c_1-b_1}\\ & = &
\dsum_{d_1=b_1}^{b_2-1}\sum_{d_2=b_2}^{b-1}
A(d_1,d_2;c_1,c_2)(1-t)^{c_1+c_2-d_1-d_2}
\\ & = & B_3(b_1,b_2,b;c_1,c_2,b) \\ & = & B_3(b_1, \dots,
b_{l-1},b; c_1, \dots, c_{l-1},b). \ea $$ If $l>3$, then by
induction $$\ba{rcl} & & \dsum_{\underline{d}\leq
\underline{c}}B_2(b_1, \dots, b_{l-1};d_1, \dots, d_{l-2},
b_{l-1})(1-t)^{\sum_{i=1}^{l-2} (c_i-d_i)} \\ & = &
\dsum_{\underline{d}\leq \underline{c}} B_3(b_1, \dots,
b_{l-1};d_1, \dots, d_{l-2},b_{l-1})(1-t)^{\sum_{i=1}^{l-2}
(c_i-d_i)}
 \\ & = & \dsum_{\underline{d}\leq
\underline{c}} \{ \sum_{d'_1=b_1}^{b_2-1}\cdots
\sum_{d'_{l-2}=b_{l-2}}^{b_{l-1}-1}A(d'_1, \dots, d'_{l-2}; d_1,
\dots, d_{l-2} )(1-t)^{\sum (d_i-d'_i)}\} (1-t)^{\sum (c_i-d_i)}
\\ & = & \dsum_{d'_1=b_1}^{b_2-1}\cdots
\sum_{d'_{l-2}=b_{l-2}}^{b_{l-1}-1}\sum_{\underline{d}\leq
\underline{c}}A(d'_1, \dots, d'_{l-2}; d_1, \dots, d_{l-2}
)(1-t)^{\sum_{i=1}^{l-2} (c_i-d'_i)} \\ & = &
\dsum_{d'_1=b_1}^{b_2-1}\cdots
\sum_{d'_{l-2}=b_{l-2}}^{b_{l-1}-1}A(d'_1, \dots, d'_{l-2},
c_{l-1}; c_1, \dots, c_{l-1} )(1-t)^{\sum_{i=1}^{l-2} (c_i-d'_i)}
\\ & = & \dsum_{d'_1=b_1}^{b_2-1}\cdots
\sum_{d'_{l-2}=b_{l-2}}^{b_{l-1}-1}\sum_{d'_{l-1}=b_{l-1}}^{b-1}
A(d'_1, \dots, d'_{l-1}; c_1, \dots, c_{l-1}
)(1-t)^{\sum_{i=1}^{l-1} (c_i-d'_i)}
\\  & = & B_3(b_1, \dots, b_{l-1},b; c_1, \dots, c_{l-1},b). \ea $$
This completes the proof of (\ref{eq2}). \par  We now assume that
$b_l<c_l$. Then $$\ba{rcl} & & B_1(b_1, \dots, b_l;c_1, \dots,
c_l)\\ & = & \dsum_{\underline{d}\in S_l} A(d_1, \dots, d_l;c_1,
\dots, c_l) \{ \ba{rcl} d_1 & \cdots & d_l
\\ b_1 & \cdots & b_l \ea \}(1-t)^{\sum (c_i-d_i)} \\ & = &
\dsum_{\underline{d}\in S_{l-1}} A(d_1, \dots, d_{l-1},c_l;c_1,
\dots, c_l) \{ \ba{cccc} d_1 & \cdots & d_{l-1} & c_l
\\ b_1 & \cdots & b_{l-1} & b_l \ea \}(1-t)^{\sum_{i=1}^{l-1}
(c_i-d_i)}. \ea $$ For every $i\leq l$, let $\sigma_i: \{b_1,
\dots, b_{l-1}\}\longrightarrow \{b_1, \dots, \hat{b}_i, \dots,
b_l\}$ be the bijection such that $\sigma_i(b_j)=b_j$ if $j<i$ and
$\sigma_i(b_j)=b_{j+1}$ if $j\geq i$. (Note that $\sigma_l$ is the
identity.) Then $$\ba{rcl} & & \{ \ba{cccc} d_1 & \cdots & d_{l-1}
& c_l
\\ b_1 & \cdots & b_{l-1} & b_l \ea \}\\ & = & \sum_{i=1}^l t
(-1)^{l+i}\{ \ba{ccc} d_1 & \cdots & d_{l-1}
\\ \sigma_i(b_1) & \cdots & \sigma_i(b_{l-1})  \ea \} \\ & = &
\sum_{i=1}^l t (-1)^{l+i} \{ \ba{cccc} d_1 & \cdots & d_{l-1} &
c_l \\ \sigma_i(b_1) & \cdots & \sigma_i(b_{l-1}) & c_l  \ea
\}.\ea $$ Therefore, $$\ba{rcl} & & B_1(b_1, \dots, b_l;c_1,
\dots, c_l)\\ & = & t\dsum_{i=1}^l(-1)^{l+i}
\dsum_{\underline{d}\in S_{l-1}} A(d_1, \dots,
d_{l-1},\underline{c})\{ \ba{cccc} d_1 & \cdots & d_{l-1} & c_l
\\ \sigma_i(b_1) & \cdots & \sigma_i(b_{l-1}) & c_l  \ea \}
(1-t)^{\sum_{i=1}^{l-1} (c_i-d_i)} \\ & = &
t\dsum_{i=1}^l(-1)^{l+i} B_1(\sigma_i(b_1), \dots,
\sigma_i(b_{l-1}),c_l;c_1, \dots, c_l). \ea $$ We have two cases
need to discuss: \par Case~1: $c_{l-1}<b_l$.
\\  In this case,
$B_1(\sigma_i(b_1), \dots, \sigma_i(b_{l-1}),c_l;c_1, \dots,
c_l)=0$ for $i<l$, as $\sigma_i(b_{l-1})=b_l>c_{l-1}$ if $i<l$.
Therefore $$B_1(b_1, \dots, b_l;c_1, \dots, c_l)=t B_1(b_1, \dots,
b_{l-1},c_l;c_1, \dots, c_l).$$ \\  On the other hand, $$\ba{rcl}
& & B_2(b_1, \dots, b_l;c_1, \dots, c_l)\\ & = &
\dsum_{\underline{d}\in S_{l-1}} A(b_1, \dots, b_l;d_1, \dots,
d_{l-1},b_l)[\ba{cccc} c_1 & \cdots & c_{l-1} & c_l
\\ d_1 & \cdots & d_{l-1} & b_l  \ea ]
(1-t)^{\sum_{i=1}^{l-1} (d_i-b_i)} \\ & = &
t\dsum_{\underline{d}\in S_{l-1}} A(b_1, \dots, b_l;d_1, \dots,
d_{l-1},b_l)[\ba{ccc} c_1 & \cdots & c_{l-1}
\\ d_1 & \cdots & d_{l-1}   \ea ]
(1-t)^{\sum_{i=1}^{l-1} (d_i-b_i)} \\ & = &
t\dsum_{\underline{d}\in S_{l-1}} A(b_1, \dots, b_{l-1},c_l;d_1,
\dots, d_{l-1},c_l)[\ba{ccc} c_1 & \cdots & c_{l-1}
\\ d_1 & \cdots & d_{l-1}  \ea ]
(1-t)^{\sum_{i=1}^{l-1} (d_i-b_i)} \\ & = &
t\dsum_{\underline{d}\in S_{l-1}} A(b_1, \dots, b_{l-1},c_l;d_1,
\dots, d_{l-1},c_l)[\ba{cccc} c_1 & \cdots & c_{l-1} & c_l
\\ d_1 & \cdots & d_{l-1} & c_l  \ea ]
(1-t)^{\sum_{i=1}^{l-1} (d_i-b_i)}\\ & = & t B_2(b_1, \dots,
b_{l-1},c_l;c_1, \dots, c_l).\ea $$ Now, we can conclude from
(\ref{eq2}) that $B_1(b_1, \dots, b_l;c_1, \dots, c_l)=B_2(b_1,
\dots, b_l;c_1, \dots, c_l)$ in the case. \par Case~2: $b_l\leq
c_{l-1}$. \\ In this case, $B_2(b_1, \dots, b_l;c_1, \dots,
c_l)=0$ by Remark~\ref{rem1}. Thus it remains to show that
$B_1(b_1, \dots, b_l;c_1, \dots, c_l)=0$. Let $\tilde \sigma_i:
\{b_1, \dots, b_l\}\longrightarrow \{b_1, \dots, \hat{b}_i, \dots,
b_l,c_l\}$ be the extension of $\sigma_i$ such that $\tilde
\sigma_i(b_j)=\sigma_i(b_j)$ for $j\leq l-1$ and $\tilde
\sigma_i(b_l)=c_l$; then by Lemma~\ref{lem6}, (\ref{eq2}) and the
fact that $c_{l-1}\geq b_l$, $$\ba{rcl}  & & B_1(b_1, \dots,
b_l;c_1, \dots, c_l)
\\ & = & t\dsum_{i=1}^l (-1)^{l+i} B_1(\sigma_i(b_1), \dots,
\sigma_i(b_{l-1}),c_l;c_1, \dots, c_l)\\ & = & t\dsum_{i=1}^l
(-1)^{l+i} B_3(\sigma_i(b_1), \dots, \sigma_i(b_{l-1}),c_l;c_1,
\dots, c_l)\\ & = &  t\dsum_{i=1}^l (-1)^{l+i} \sum_{d_1=\tilde
\sigma_i(b_1)}^{\tilde \sigma_i(b_2)-1}\cdots \sum_{d_{l-1}=
\tilde \sigma_i(b_{l-1})}^{\tilde \sigma_i(b_l)-1} A(d_1, \dots,
d_{l-1};c_1, \dots, c_{l-1}) (1-t)^{\sum (c_i-d_i)} \\ & = &
t\dsum_{d_1=b_1}^{b_2-1}\cdots \sum_{d_{l-1}= b_{l-1}}^{b_l-1}
A(d_1, \dots, d_{l-1};c_1, \dots, c_{l-1}) (1-t)^{\sum
(c_i-d_i)}\\ & = & 0. \ea $$ The proof of the theorem is now
complete.

\section{Main theory}
The way to prove our main theorem is to use induction on $r$.
Therefore we shall first discuss the case $r=2$.
\begin{prop}
Let $Y$ be a one-sided ladder-shaped region of the plane with the
partial order given in the beginning. Let $P_i=(i, n)$ and $Q=(m,
b_i)$, $i=1,2$ be points of $Y$. Let $W(n,m; b_1, b_2)$ be the
polynomial $\sum_{k\geq 0} w_k({\bf P}, {\bf Q})t^k$ and let
$\tilde W(n,m; b_1, b_2)=\det[W(P_i, Q_j)]_{i,j=1,2}$. Then
$$t^{-1}\tilde W(n,m; b_1, b_2)= \sum_{(c_1,c_2)\in S_2}
A(b_1,b_2;c_1,c_2)(1-t)^{c_1+c_2-b_1-b_2}W(n,m; c_1,c_2).$$ In
particular, if $b_1=1$ and $b_2=2$, then $$t^{-1}\tilde W(n,m;
1,2)= W(n,m; 1,2).$$ \label{prop1}
\end{prop}
\begin{proof}
We prove by induction on $m$. If $m=2$, then $$\ba{rcl} & &
\sum_{(c_1,c_2)\in S_2}
A(b_1,b_2;c_1,c_2)(1-t)^{c_1+c_2-b_1-b_2}W(n,2; c_1,c_2)
\\ & = &
\sum_{c_1=b_1}^{b_2-1}(1-t)^{c_1-b_1}[b_2-c_1-(b_2-c_1-1)(1-t)]
\\ & = & b_2-b_1 \\& = & t^{-1} \tilde W(n,2;b_1,b_2). \ea $$ Assume
that $m>2$. Then by Lemma~\ref{lem2} and induction $$\ba{rcl} & &
t^{-1}\tilde W(n,m; b_1,b_2)\\  & = & t^{-1} \dsum_{(c_1,c_2)\in
S_2,(m,c_2)\in Y} \{ \ba{cc} c_1 & c_2 \\ b_1 & b_2 \ea \} \tilde
W(n,m-1;c_1,c_2) \\ & = & \dsum_{(c_1,c_2)\in S_2,(m,c_2)\in Y} \{
\ba{cc} c_1 & c_2 \\ b_1 & b_2 \ea \} \sum_{(d_1,d_2)\in
S_2,(m,d_2)\in Y} A(c_1,c_2;d_1,d_2) W(n,m-1;d_1,d_2) \\ & = &
\dsum_{(d_1,d_2)\in S_2,(m,d_2)\in Y}[ \sum_{(c_1,c_2)\in
S_2}A(c_1,c_2;d_1,d_2) \{ \ba{cc} c_1 & c_2 \\ b_1 & b_2 \ea \}]
W(n,m-1;d_1,d_2) \\ & = & \dsum_{(d_1,d_2)\in S_2,(m,d_2)\in
Y}B_1(b_1,b_2;d_1,d_2)  W(n,m-1;d_1,d_2). \ea $$ Furthermore, by
Lemma~\ref{lem2} $$\ba{rcl} &  & \dsum_{(c_1,c_2)\in S_2}
A(b_1,b_2;c_1,c_2)(1-t)^{c_1+c_2-b_1-b_2}W(n,m; c_1,c_2) \\ & = &
\dsum_{(c_1,c_2)\in S_2}
A(b_1,b_2;c_1,c_2)(1-t)^{c_1+c_2-b_1-b_2}\sum_{(d_1,d_2)\in
S_2,(m,d_2)\in Y} [\ba{cc} d_1 & d_2 \\ c_1 & c_2  \ea ]
W(n,m-1;d_1,d_2)\\ & = & \dsum_{(d_1,d_2)\in S_2,(m,d_2)\in
Y}\{\sum_{(c_1,c_2)\in S_2}
A(b_1,b_2;c_1,c_2)(1-t)^{c_1+c_2-b_1-b_2}[\ba{cc} d_1 & d_2 \\ c_1
& c_2  \ea ]  \} W(n,m-1;d_1,d_2) \\ & = & \dsum_{(d_1,d_2)\in
S_2,(m,d_2)\in Y}B_2(b_1,b_2;d_1,d_2)  W(n,m-1;d_1,d_2). \ea $$
Since $B_1(b_1,b_2;d_1,d_2)=B_2(b_1,b_2;d_1,d_2)$, the assertion
follows.
\end{proof}

\begin{theo}
Let $Y$ be a one-sided ladder-shaped region of the plane (see
Figure~1) with the partial order given in the beginning. Let
$P_i=(i, n)$ and $Q=(m, b_i)$, $i=1, \dots, r$ be points of $Y$.
Let $W(n,m; b_1, \dots, b_r)$ be the polynomial $\sum_{k\geq 0}
w_k({\bf P}, {\bf Q})t^k$ and let $\tilde W(n,m; b_1, \dots,
b_r)=\det[W(P_i, Q_j)]_{i,j=1, \dots, r}$. Then $$t^{-{r\choose
2}}\tilde W(n,m; b_1, \dots, b_r)= \sum_{\underline{c}\in S_r}
A(\underline{b}; \underline{c})(1-t)^{\sum (c_i-b_i)}W(n,m; c_1,
\dots, c_r),$$ where $\underline{b}=(b_1, \dots, b_r)$ and
$\underline{c}=(c_1, \dots, c_r)$. In particular, if
$b_i=i~\forall i$, then $$t^{-{r\choose 2}}\tilde W(n,m; 1, \dots,
r)= W(n,m; 1, \dots, r).$$ \label{theo4}
\end{theo}

\begin{proof}
We prove the theorem by induction on $r$ and $m$. If $r=2$, then
it is the content of Proposition~\ref{prop1}. Therefore we may
assume that $r>2$. \\ Assume for the moment that $m=r$. Notice
that by Lemma~\ref{lem2} $$\ba{rcl} & & \dsum_{\underline{c}\in
S_r} A(\underline{b}; \underline{c})(1-t)^{\sum (c_i-b_i)}W(n,m;
c_1, \dots, c_r)\\ & = & \dsum_{\underline{c}\in S_{r-1},c_r=b_r}
A(\underline{b}; \underline{c})(1-t)^{\sum (c_i-b_i)}W(b_r-1,r;
c_1, \dots, c_{r-1}) \\ & = & \dsum_{\underline{c}\in
S_{r-1},c_r=b_r} A(\underline{b}; \underline{c})(1-t)^{\sum
(c_i-b_i)}\{\sum_{\underline{d}\in S_{r-1},d_{r-1}<b_r} [\ba{ccc}
d_1 & \cdots & d_{r-1} \\ c_1 & \cdots & c_{r-1} \ea ]
W(b_r-1,r-1; \underline{d})\}\\ & = & \dsum_{\underline{d}\in
S_{r-1},d_{r-1}<b_r}\{\sum_{\underline{c}\in S_{r-1},c_r=b_r}
A(\underline{b}; \underline{c})[\ba{ccc} d_1 & \cdots & d_{r-1}
\\ c_1 & \cdots & c_{r-1} \ea ](1-t)^{\sum (c_i-b_i)}\}
W(b_r-1,r-1; \underline{d})\} \\ & = & \dsum_{\underline{d}\in
S_{r-1},d_{r-1}<b_r} B_2(b_1, \dots, b_r;d_1, \dots, d_{r-1},b_r
). \ea $$ Furthermore, by Corollary~\ref{coll1} and induction,
$$\ba{rcl} & & t^{-{r\choose 2}}\tilde W(n,m; b_1, \dots, b_r) \\
& = & t^{-{r-1\choose 2}}\dsum_{c_1=b_1}^{b_2-1} \cdots
\sum_{c_{r-1}=b_{r-1}}^{b_r-1} \tilde W(b_r-1,r-1; c_1, \dots,
c_{r-1})\\ & = & \dsum_{c_1=b_1}^{b_2-1} \cdots
\sum_{c_{r-1}=b_{r-1}}^{b_r-1}\sum_{\underline{d}\in
S_{r-1}}A(\underline{c},\underline{d})(1-t)^{\sum
(d_i-c_i)}W(b_r-1,r-1;d_1, \dots, d_{r-1}) \\ & = &
\dsum_{\underline{d}\in
S_{r-1},d_{r-1}<b_r}\{\sum_{c_1=b_1}^{b_2-1} \cdots
\sum_{c_{r-1}=b_{r-1}}^{b_r-1}A(\underline{c},\underline{d})(1-t)^{\sum
(d_i-c_i)} \}W(b_r-1,r-1;d_1, \dots, d_{r-1}) \\ & = &
\dsum_{\underline{d}\in S_{r-1},d_{r-1}<b_r} B_3(b_1, \dots,
b_r;d_1, \dots, d_{r-1},b_r ).  \ea $$ Thus the theorem holds for
$m=r$ by Theorem~\ref{theo3}. \\ We now assume that $m>r$. By
Lemma~\ref{lem2} and induction, $$\ba{rcl} & &  t^{-{r\choose
2}}\tilde W(n,m; b_1,\dots, b_r)\\ & = & t^{-{r\choose 2}}
\dsum_{\underline{c}\in S_r,(m,c_r)\in Y} \{ \ba{ccc} c_1 & \cdots
& c_r
\\ b_1 & \cdots & b_r \ea \} \tilde W(n,m-1;c_1,\dots, c_r) \\ & = &
\dsum_{\underline{c}\in S_r,(m,c_r)\in Y} \{ \ba{ccc} c_1 & \cdots
& c_r
\\ b_1 & \cdots & b_r \ea \} \sum_{\underline{d}\in S_r,(m,d_r)\in Y}
A(\underline{c};\underline{d})(1-t)^{\sum (d_i-c_i)}
W(n,m-1;d_1,\dots, d_r) \\ & = & \dsum_{\underline{d}\in
S_r,(m,d_r)\in Y} [\sum_{\underline{c}\in S_r,(m,c_r)\in Y}
A(\underline{c};\underline{d})(1-t)^{\sum (d_i-c_i)} \{ \ba{ccc}
c_1 & \cdots & c_r
\\ b_1 & \cdots & b_r \ea \}] W(n,m-1;d_1,\dots, d_r)
\\ & = & \dsum_{\underline{d}\in S_r,(m,d_r)\in
Y}B_1(\underline{b};\underline{d}) W(n,m-1;d_1,\dots, d_r).\ea $$
Also, by Lemma~\ref{lem2}  $$\ba{rcl} & & \dsum_{\underline{c}\in
S_r} A(\underline{b}; \underline{c})(1-t)^{\sum (c_i-b_i)}W(n,m;
c_1, \dots, c_r) \\ & = & \dsum_{\underline{c}\in S_r}
A(\underline{b}; \underline{c})(1-t)^{\sum (c_i-b_i)}
\sum_{\underline{d}\in S_r,(m,d_r)\in Y} [ \ba{ccc} d_1 & \cdots &
d_r
\\ c_1 & \cdots & c_r \ea ]
W(n,m-1;d_1,\dots, d_r)\\ & = & \dsum_{\underline{d}\in
S_r,(m,d_r)\in Y}\{\sum_{\underline{c}\in S_r} A(\underline{b};
\underline{c})(1-t)^{\sum (c_i-b_i)} [ \ba{ccc} d_1 & \cdots & d_r
\\ c_1 & \cdots & c_r \ea ]  \} W(n,m-1;d_1,\dots, d_r) \\ & = &
\dsum_{\underline{d}\in S_r,(m,d_r)\in Y}B_2(\underline{b};
\underline{d})W(n,m-1;d_1,\dots, d_r).\ea $$ Since
$B_1(\underline{b}; \underline{d})=B_2(\underline{b};
\underline{d})$, the assertion follows.
\end{proof}

\end{document}